\documentclass[11pt]{amsart}
\usepackage{amsrefs}
\usepackage[margin=0.8in]{geometry}
\usepackage{amssymb,color}
\usepackage{amsmath}
\usepackage{graphicx}
\usepackage{hyperref,ulem}
\usepackage{amssymb}
\usepackage{MnSymbol}
\usepackage{enumerate}
\usepackage[english]{babel}
\newtheorem{theorem}{Theorem}[section]
\newtheorem{lemma}[theorem]{Lemma}
\newtheorem{proposition}[theorem]{Proposition}

\newtheorem{definition}[theorem]{Definition}
\newtheorem{remark}{Remark}
\newcommand{\R}{\mathbb{R}}

\newcommand{\id}{\mathrm{id}}
\newcommand{\I}{\mathrm{i}}
\newcommand{\e}{\operatorname{e}}

\newcommand{\SF}{ shear-free }

\title{On a criterion for local embeddability of 3-dimensional CR-structures}
\author{Masoud Ganji}
\email[M. Ganji]{mganjia2@une.edu.au}
\author{Gerd Schmalz}
\email[G. Schmalz]{schmalz@une.edu.au}
\address{School of Science and Technology, University of New England, Armidale NSW 2351, Australia}
\begin{document}
\maketitle
\begin{abstract} We introduce a CR-invariant class of Lorentzian metrics on a circle bundle over a 3-dimensional CR-structure, which we call quasi-Fefferman metrics. These metrics generalise the Fefferman metric but allow for more control of the Ricci curvature. Our main result is a criterion for embaddability of 3-dimensional CR-structures in terms of the Ricci curvature of the quasi-Fefferman metrics in the spirit of the results by Lewandowski et al. in \cite{MR1078890} and also Hill et al. in \cite{MR2492229}.   
\end{abstract}
\section{Introduction}
A 3-dimensional abstract strictly pseudoconvex CR-structure is a 3-dimensional manifold $M$ with a contact distribution $H$ and a smooth field of endomorphisms $J_x\colon H_x\to H_x$ such that $J_x^2=-\id$. Such structures naturally occur on strictly pseudoconvex real hypersurfaces $M$ in $\mathbb C^2$ where 
$$H_x=T_xM\cap J_xT_xM$$
and $J_x$ is the canonical complex structure on $T_x\mathbb C^2$. The (local) embeddability (or realisability) problem asks if there exists a (local) embedding $\iota: M \to \mathbb C^2$ of the abstract CR-manifold $M$ such that the CR-structures on $\iota M$ induced by $\iota$ and by $\mathbb C^2$ coincide. This problem is equivalent to finding two functionally independent CR-functions, that is solutions to the complex linear PDE
\begin{equation}\label{hle}
\bar{\partial} f=0,
\end{equation}
where $\bar{\partial}= X+\I JX$ for some non-vanishing section $X$ of $H$.

While equation \eqref{hle} has sufficiently many solutions for real-analytic CR-structures (see, e.g., \cite{MR1067341}) to make them embeddable it is well known that for almost all non-analytic CR-structures only constant CR-functions exist \cite{MR656622, MR0450756, MR0410042}.

It has been known to physicists since the 1960s (see e.g. \cite{MR0250641}) that the existence of a non-constant CR-function on $M$ is equivalent to the existence of a Lorentzian metric $g$ on $\mathcal M=M\times \mathbb R$ that satisfies two properties
\begin{itemize}
\item[a.)]  The fundamental vector field $k=\partial_r$ of the line bundle $M\times \mathbb R$ is null and shear-free with respect to $g$, that is  $g(k,k)=0$ and  
\begin{equation}\label{sfv}
\mathcal L_k g= \rho g + \theta\vee \psi,
\end{equation}
where $\mathcal L$ is the Lie derivative,  $\rho$ is a function on $\mathcal M$, $\theta=g(k,\cdot)$ and $\psi$ is some 1-form on $\mathcal M$.
\item[b.)] The complexified Ricci tensor of $g$ vanishes on $\alpha$-planes, i.e. $g$ is partially Ricci flat. (See \eqref{plane} below for the precise definitions. )
\end{itemize}
It turns out that the existence of non-constant solutions to \eqref{hle} is not sufficient to guarantee embeddability of the CR-structure \cite{MR1030851}. A remarkable theorem by Jacobowitz \cite{MR876018} provides a criterion for embeddability of CR-structures in terms of the canonical bundle (see Definition \eqref{defcan}). Below we reformulate Jacobowitz's more general result for 3-dimensional CR-structures. 

\begin{theorem}[Jacobowitz , 1989]\cite{MR876018}\label{cls} 
Let  $(M, H, J)$ be a CR-structure. Suppose that near some point $p\in M$ the CR-structure has a non-constant CR-function. If the canonical bundle  associated with the CR-structure has a non-vanishing closed  section then the CR-structure is embeddable near $p$.
\end{theorem}

Notice that the converse statement of Theorem \ref{cls} is also true: For an embeddable CR-structure we have two functionally independent, hence non-constant, coordinate CR-functions $\zeta,\eta$. The 2-form $d\zeta\wedge d\eta$ gives a non-vanishing closed section of the canonical bundle.  

Lewandowski, Nurowski and Tafel \cite{MR1078890} and also Hill, Lewandowski and Nurowski \cite{MR2492229} prove a series of embeddability results in terms of shear-free congruences of Lorentzian spaces. In particular, they show that a CR-structure is embeddable if and only if it admits a lift to a partially Ricci flat Lorentz space $(\mathcal M, g)$ with a \SF congruence $k$ as above, and a Maxwell field aligned with $k$. This is closely related to Jacobowitz's theorem \ref{cls}, since the existence of such Maxwell field is equivalent to the existence of a non-vanishing closed section of the canonical bundle from Theorem \eqref{cls}.

In this paper we introduce a family of Lorentzian metrics on a circle bundle over the CR-manifold $M$. This family is more general than the conformal class of Fefferman metrics but more special than the family of Lorentzian metrics that admit a \SF congruence. We call these metrics \textit{ quasi-Fefferman metrics}. Our main result is the following embeddability criterion.  

\begin{theorem}\label{main}
A strictly pseudoconvex CR-structure $(M,H,J)$ is locally embeddable if and only if it admits a quasi-Fefferman metric for which the complexified Ricci curvature vanishes on the $\alpha$-planes. 
\end{theorem}

Notice that the special choice of the metrics allows us to drop the assumption of the existence of a non-vanishing closed section of the canonical bundle.  In order to give a self-contained proof we recapitulate techniques used in \cite{MR2492229} and references therein.

\section{CR-structures and the Fefferman metric}
Let $(M,H,J)$ be 3-dimensional CR-structure. We assume that $(M,H,J)$ is strictly pseudoconvex, i.e. for any (local) non-vanishing section $X$ of $H$, $[X,JX]\not\in H$. For any choice of $X$ we have an adapted complex frame $\partial,\bar{\partial},\partial_0$, where
$$\partial= X-\I JX,\qquad \partial_0=\I[\partial,\bar{\partial}]=-2[X,JX].$$
The complex vector field $\partial= X-\I JX$ spans the $+\I$-eigendistribution $H^{1,0}$ of $J$ in $H\otimes\mathbb C$.
We denote the 
corresponding dual coframe by $(\mu,\bar{\mu},\lambda)$. Strict pseudoconvexity of $M$ translates to
$$d\lambda \wedge \lambda \neq 0.$$ 
Our choice implies
\begin{align}\label{cofcon}
d\lambda&=\mathrm{i} \mu \wedge \bar{\mu}+c\mu \wedge \lambda+\bar{c}\bar{\mu} \wedge \lambda\\\nonumber
d\mu&=\alpha\mu \wedge \lambda +\beta \bar{\mu} \wedge \lambda,
\end{align}
where $c, \alpha, \beta$ are complex-valued functions on $M$. Any other distinguished frame $(\partial',\bar{\partial'},\partial_0')$ and coframe $(\mu',\bar{\mu'},\lambda')$ express through the original frame and coframe by
\begin{align}\label{trans1}
\partial'&=\e^{-\tau-\I \theta}\partial,& 
\qquad \partial_0'&=\e^{-2\tau}\big(\partial_0-h\partial-\bar{h}\bar{\partial}\big) \\\nonumber
\mu'&=\e^{\tau+\I \theta}(\mu+h\lambda), &\qquad \lambda'&=\e^{2\tau}\lambda
\end{align}
where $\tau$ and $\theta$ are real-valued functions and
\begin{align}\label{equa11}
&h=-\I \bar{\partial}(\tau+\I \theta), 
\quad \alpha'=\e^{2\tau}
\left(\alpha-\partial_0 (\tau+\I \theta)+h\partial (\tau+\I \theta)+\partial h+hc \right), \quad c'=\e^{-\tau-\I \theta}\left(c-2\I\bar{h}+\partial (\tau+\I \theta)\right).
\end{align}

 Recall that the Fefferman metric is a conformal class of Lorentzian metrics defined on the circle bundle $\mathcal M= H^{1,0}/{\mathbb R^+}$. Using a distinguished coframe $(\mu, \lambda)$ and the trivialisation $\mathcal M \ni m|_p= \e^{-t-\I \rho} \partial|_p\mapsto (p,\rho)$ where $\rho\in [0, 2 \pi)$, 
a representative of the Fefferman metric is defined by the simple formula
\begin{align}\label{feff}
g_F=  
\mu  \bar{\mu}   + \lambda\left(\frac{2}{3} d\rho -\frac{\I}{3}c  \mu +\frac{\I}{3}\bar{c} \bar{\mu} - \left(\frac{\partial \bar{c}+\bar{\partial}c}{12}-\frac{\I(\alpha-\bar{\alpha})}{4}\right)\lambda \right)
\end{align}
where we kept the notations $\mu,\bar{\mu},\lambda$ for their pull-backs under the circle-bundle projection, and $\alpha$ and $c$ are as above (see \cite{MR943984}). The CR-invariance of the Fefferman metric means that a change of the distinguished coframe $(\mu,\bar{\mu},\lambda)$ causes only a conformal change of $g_F$ by the factor $\e^{2\tau}$.

Denote by $\tilde{\mathcal M}$ the natural lift of $\mathcal M$ to a line bundle. It will be convenient in the computations below to rescale the coordinate $\rho$ on $\tilde{\mathcal M}$ to $r=\frac{2\rho}{3}$. Then the change of the coframe $(\lambda,\mu)$ induces the change 
\begin{equation}\label{rcha}
r'=r-\frac{2}{3}\theta
\end{equation}

of the trivialisation of $\tilde{\mathcal M}$. Denote the quotient bundle of the (rescaled) line bundle $\tilde{\mathcal M}$ mod $2\pi$ by $\mathcal M^{\frac32}$. Since the Fefferman metric is invariant with respect to the principle $\mathbb R$-bundle action it projects to any $S^1$-bundle with arbitrary period. In particular, it is well defined on $\mathcal M^{\frac32}$.  

It is known that the Weyl tensor of the Fefferman metric has rank (at most) one \cite{MR943984}. However, in general, the Ricci curvature cannot be controlled and, except for very special cases (see \cite{MR2289882, MR2463976}), one cannot find (partially) Ricci flat or Einstein representatives of the conformal Fefferman metric. 

\section{Quasi-Fefferman metrics}
We generalise the Fefferman metric by introducing functional parameters $x$  and $H$ in addition to the conformal factor. This allows us to impose conditions on the Ricci curvature that cannot be satisfied by the Fefferman metric. We will see later that it is more natural to define the quasi-Fefferman metrics on the circle bundle $\mathcal M^{\frac32}$ rather than the Fefferman bundle $\mathcal M$.
\begin{definition}
Let $(M,H,J)$ be a CR-structure and  $\mathcal  M^{\frac32}$ as above. For any choice of a distinguished coframe $(\mu,\lambda)$ and the induced trivialisation of $\mathcal M^{\frac32}$ we define the family of quasi-Fefferman metrics on $\mathcal  M^{\frac32}$ by 
\begin{equation}\label{metric1}
g=2P^2\left[\mu \bar{\mu}+\lambda\left( dr+W\mu+\overline{W}\bar{\mu}+H\lambda \right)\right]
\end{equation}
 where  
$W=\I x\e^{-\I 	 r}-\frac{\I}{3} c.$
Here $P\neq 0, H$ are real-valued functions on $\mathcal M^{\frac32}$ and  $x$ is a complex-valued function on $M$.
\end{definition}

It is an important feature of the family of \SF metrics \eqref{metric1} and of the conformal family of Fefferman metrics that they are CR invariant, i.e. they do not depend on the choice of the pair $(\mu,\lambda)$. 
We show that this is also true for the family of quasi-Fefferman metrics.
\begin{theorem}
The family of quasi-Fefferman metrics is CR invariant. 
\end{theorem}
{\bf{Proof.}} Under the frame change \eqref{trans1} and the induced change of the trivialisation \eqref{rcha} the quasi-Fefferman metric changes as follows.
\begin{align*}
g'&=2P'^2\left[\mu'\bar{\mu'}+ \lambda'\left( 	dr'
+ (\I x' \e^{-\I r'}-\frac{\I}{3}c')  \mu'
 +(-\I \bar{x'} \e^{\I r'}+\frac{\I}{3} \bar{c'})\bar{\mu'} 
 +H'\lambda' \right)\right].\\
g'&=2|f|^2P'^2\bigg[\mu\bar{\mu}+\bar{h}\mu \lambda+h \bar{\mu}\lambda+|h|^2 \lambda^2+ \lambda\bigg( dr-\frac{2}{3}d\theta
+ f\left (\I x' \e^{\I\frac{2}{3} \theta}\e^{-\I r}-\frac{\I}{3f}(c-2\I\bar{h}+\partial \log f)\right)  \mu\\
 &+\bar{f} \left (-\I \bar{x'} \e^{-\I\frac{2}{3} \theta}\e^{\I r}+\frac{\I}{3\bar{f}} (\bar{c}+2\I h+\bar{\partial} \log \bar{f})\right )\bar{\mu}
 +\left (fh(\I x' \e^{\I\frac{2}{3} \theta}\e^{-\I r}-\frac{\I}{3}c')
 +\bar{f}\bar{h}(-\I \bar{x'} \e^{-\I\frac{2}{3} \theta}\e^{\I r}+\frac{\I}{3} \bar{c'})+|f|^2H'\right)\lambda \bigg)\bigg].\\
&=2P^2\bigg[\mu\bar{\mu}+ \lambda\bigg( dr
+ (\I x \e^{-\I r}-\frac{\I}{3}c) \mu
+(\frac{1}{3}\bar{h}-\frac{\I}{3}\partial \log f-\frac{2}{3}\partial \theta)  \mu
 +  (-\I \bar{x} \e^{\I r}+\frac{\I}{3} \bar{c})\bar{\mu}+(\frac{1}{3} h+\bar{\partial} \log \bar{f}-\frac{2}{3}\bar{\partial} \theta )\bar{\mu}
+H\lambda \bigg)\bigg].\\
 &=2P^2\bigg[\mu\bar{\mu}+ \lambda\bigg( dr
+ (\I x \e^{-\I r}-\frac{\I}{3}c) \mu
+  (-\I \bar{x} \e^{\I r}+\frac{\I}{3} \bar{c})\bar{\mu}+H 
\lambda \bigg)\bigg].
%
%
%
\end{align*}
where 
\begin{align}\label{qfch}
&f=\e^{\tau+\I \theta}, \quad P=\e^{\tau}P', \quad x=\e^{\tau+\I \frac{5}{3}\theta}	x', 
\\\nonumber
  &H=\e^{2\tau}H'+|h|^2+\e^{\tau+\I \theta}h(\I x' \e^{-\I r'}-\frac{\I}{3}c')
 +\e^{\tau-\I \theta}\bar{h}(-\I \bar{x'} \e^{\I r'}+\frac{\I}{3} \bar{c'})-\frac{2}{3}\partial_0 \theta. \qquad \Box 
 \end{align}
\section{Lorentzian geometry and $\alpha$-planes}
Let $(\mathcal M, g)$ be a 4-dimensional Lorentzian manifold  equipped with a foliation into integral curves of a non-vanishing null vector field $k$. We have the following canonical objects
\begin{itemize}
\item[(i)] the 1-form $\theta=g(k,\cdot)$
\item[(ii)] the distribution $k^\perp=\{X\in \Gamma(T\mathcal M)\colon g(X,k)=0\}$
\item[(iii)] the distribution of \textit {screen spaces} $S:=k^\perp/k$.
\end{itemize}

\begin{proposition}
On each screen space $S_x$ there are two canonical complex structures $J_x$ and $-J_x$.
\end{proposition}

{\bf Proof.} The restriction of $g$ to $k^\perp$ is a degenerate metric with kernel $k$ and induces a Euclidean metric on $S_x$. Since $S_x$ is 2-dimensional the Euclidean metric induces complex structures of rotation by $\frac{\pi}{2}$ in either orientation. \hfill $\Box$\medskip

Choose one of the two complex structures on $S$. Then $\mathbb C \otimes S$ splits into its eigenspaces $S^{1,0} \oplus S^{0,1}$.
Let 
$$\pi\colon \mathbb C \otimes k^\perp \to \mathbb C \otimes S$$
be the canonical projection map. 
The subspaces $K^{1,0}, K^{0,1}$ of $\mathbb C \otimes k^\perp$  defined by 
\begin{equation}\label{plane}
K^{1,0}=\pi^{-1}S^{1,0}, \quad K^{0,1}=\pi^{-1}S^{0,1}
\end{equation}
are called $\alpha$-planes  and $\beta$-planes, respectively. Notice that changing the orientation used in the definition of $J$ results in interchanging the $\alpha$-planes and the $\beta$-planes.
 
Clearly, 
$$K^{1,0} \cap K^{0,1}=k, \quad K^{1,0} + K^{0,1}=\mathbb C \otimes k^\perp.$$

\begin{definition}
We say that the complexified Ricci tensor of $g$ vanishes on the $\alpha$-planes $K^{1,0}$ if $Ric|_{K^{1,0}}=0$, i.e.,
$$\text{Ric}(X_1, X_2)=0 \quad \forall X_1, X_2 \in K^{1,0}. $$
\end{definition} 
 
Notice that vanishing of the complexified Ricci tensor on $\alpha$-planes is equivalent to its vanishing on $\beta$-planes. Hence the definition above does not depend on the choice of $J$. 
 
\begin{definition}
Let $\mathcal M$ be a 4-dimensional manifold equipped with a Lorentzian metric $g$ and a non-vanishing null vector field $k$. A complex frame $(e_1, e_2, \ell, k)$ is called adapted to $(g,k)$ if $e_1$ is a section of $\alpha$-planes, $e_2=\bar{e}_1$ (and, hence, is a section of $\beta$-planes), and
$$g(e_1,e_2)=1, \quad g(\ell,\ell)=0, \quad g(\ell,k)=1, \quad g(\ell,e_1)=\overline{g(\ell,g_2)}=0.$$  
\end{definition} 

\begin{proposition}
A 4-dimensional Lorentzian manifold $(\mathcal{M},g)$ with a non-vanishing null vector field $k$ possesses (locally) a complex adapted frame.  
\end{proposition}

{\bf Proof.} Let $\ell$ be a null vector field such that $g(\ell,k)=1$. Choose a unit vector field $\varepsilon_1\in k^\perp \cap \ell^\perp$. Choose $\varepsilon_2\in k^\perp$ such that $\pi{\varepsilon_2}=J\pi(\varepsilon_1)$ and $g(\varepsilon_2,\ell)=0$. 
Now, set 
$$e_1= \frac{1}{\sqrt{2}}(\varepsilon_1 - \I \varepsilon_2), \quad e_2= \frac{1}{\sqrt{2}}(\varepsilon_1 + \I \varepsilon_2). \qquad \Box$$
 \medskip
 
It follows that the $\alpha$-planes are spanned by $(e_1,k)$ and the $\beta$-planes are spanned by $(e_2,k)$.  
 
Now, vanishing of the Ricci curvature on $\alpha$-planes is equivalent to
\begin{itemize}
\item[(i)] $R_{11}=Ric(e_1,e_1)=\overline{Ric (e_2,e_2)}=0$
\item[(ii)] $R_{14}=Ric(e_1,k)=\overline{Ric (e_2,k)}=0$
\item[(iii)] $R_{44}=Ric(k,k)=0$.
\end{itemize}
   
For the dual complex coframe $(\theta^1, \theta^2, \theta^3, \theta^4)$ to an adapted complex frame 
 $\bar{\theta^1}=\theta^2$ vanishes on $\alpha$-planes, $\theta^3=g(k,\cdot)$ and 
  $$g=2(\theta^1\theta^2+\theta^3\theta^4).$$
The Gram matrix for both $g$ and its dual $g^{-1}$ with respect to an adapted frame and coframe is:
\begin{equation}\label{gram}
\begin{pmatrix} 0& 1& 0&0\\ 1& 0 & 0 & 0\\ 0& 0 & 0 & 1\\ 0 & 0 & 1 & 0 \end{pmatrix}.
\end{equation}
  
A direct computation shows that the shearfreeness condition \eqref{sfv} is equivalent to
\begin{equation}\label{sfv1}
 d\theta^3\wedge \theta^1 \wedge \theta^3=0, \quad d\theta^1\wedge \theta^1 \wedge \theta^3=0.
 \end{equation}
See, e.g.,  \cite{MR2492229} for details.

Below we cite a version of the celebrated Goldberg-Sachs theorem \cite{MR0156679, MR2786175, MR2492229}, which is a useful tool for computing certain components of the Weyl tensor of the Lorentzian metric $g$:
$$C_{ijkl}=R_{ijkl}+\frac{1}{6}R\big (g_{ik}g_{lj}-g_{il}g_{kj}\big )+\frac{1}{2}\big(g_{il}R_{kj}-g_{ik}R_{lj}+g_{jk}R_{li}-g_{jl}R_{ki}\big)$$
where $R_{ijkl}$ is the Riemann curvature, $R_{kj}$ is the Ricci curvature and $R$ is the scalar curvature. The following quantities are called Weyl scalars:
$$\Psi_0= C(k, e_1, k, e_1)=
C_{4141},
 \qquad \Psi_1=
C(k, l, k, e_1)=
C_{4341}
.$$

\begin{theorem}[Goldberg-Sachs theorem, \cite{MR2786175, MR2492229} ]\label{goldsachs}
 Suppose that a 4-dimensional manifold $\mathcal M$ is equipped with a Lorentzian metric $g$
  and a shearfree null vector field $k$, as above.
 Also assume that the complexified Ricci curvature of $g$ vanishes on $\alpha$-planes, i.e. $R_{11}= R_{14}= R_{44}= 0$ with respect to an adapted coframe $(\theta^1, \theta^2,  \theta^3, \theta^4)$. Then the Weyl scalars $\Psi_0=\Psi_1=0$.
 \end{theorem}

The next lemma plays a crucial role in finding a CR function. It reduces the existence of a non-vanishing CR-function to the existence of a certain complex-valued 1-form. Such 1-form can be obtained from the Levi-Civita connection form of the associated Lorentzian metric, if the Ricci curvature vanishes on $\alpha$-planes.  The proof is based on Frobenius's theorem. For a detailed proof see e.g.  \cite{MR2492229}.
\begin{lemma}\label{firstcr}\cite{MR2492229}
Let $\varphi$ be a smooth complex valued 1-form defined locally in $\R^n$, $n \geq 3$, such that $\varphi \wedge \bar{\varphi} \neq 0$. Then
$$d \varphi \wedge \varphi =0$$ 
 if and only if there exist  smooth complex functions $\zeta$ and $h$  such that 
 $$\varphi=hd\zeta, \qquad d \zeta \wedge d\bar{\zeta} \neq 0.$$ 
\end{lemma}

\section{Quasi- Fefferman metrics and the embedding of 3-dimensional CR-structures}
Before we prove the main theorem of this paper we collect some ingredients. First we compute the Levi-Civita connection 1-forms, $\Gamma^i_j$ with respect to an adapted frame.
Notice that 
due to \eqref{cofcon}  we have
$$[\partial, \partial_0]=-\alpha \partial-\bar{\beta}\bar{\partial}-c\partial_0, \qquad [\bar{\partial}, \partial_0]=-\beta \partial-\bar{\alpha}\bar{\partial}-\bar{c}\partial_0. $$
Now in terms of the coframe $(\theta^1, \theta^2, \theta^3, \theta^4)$ with
\begin{align}\label{coframe}
 {\theta}^1=P\mu, \qquad {\theta}^2=P\bar{\mu}, \qquad {\theta}^3=P\lambda, \qquad {\theta}^4=P\big (dr+W\mu+\overline{W}\bar{\mu}+H\lambda \big).
 \end{align}
the metric \eqref{metric1} becomes
\begin{align}\label{metric2}
g=2\theta^1\theta^2+2\theta^3\theta^4.
\end{align}
The dual frame $(e_1, e_2, e_3, e_4)$ to $(\theta^1, \theta^2, \theta^3, \theta^4)$ takes the form
\begin{equation}\label{frame}
e_1=\frac{1}{P}(\partial-W{\partial}_r), \qquad e_2=\frac{1}{P}(\bar{\partial}-\overline{W}{\partial}_r), \qquad e_3=\frac{1}{P}(\partial_0-H {\partial}_r), \qquad e_4=\frac{1}{P}{\partial}_r.
\end{equation}
The commutators of the frame field \eqref{frame}  evaluate to
\begin{align*}
&[e_1, e_2]=\left(\frac{\bar{\partial}P}{P^2}-\overline{W}\frac{P_r}{P^2}\right)e_1+\left(-\frac{\partial P}{P^2}+W \frac{P_r}{P^2}\right)e_2-\frac{\I}{P}e_3
+\left(-\frac{\I H}{P}+ W_2-\overline{W}_1\right)e_4\\
&[e_1, e_3]=\left(\frac{\partial_0 P}{P^2}-H \frac{P_r}{P^2}	-\frac{\alpha}{P}\right)e_1-\frac{\bar{\beta}}{P}e_2+\left(-\frac{\partial P}{P^2}+W \frac{P_r}{P^2}	-\frac{c}{P}\right)e_3+
 \left(-\frac{cH}{P}+ W_3-H_1-\frac{\alpha W}{P}-\frac{\bar{\beta}\overline{W}}{P} \right)e_4\\
&[e_1, e_4]=\frac{P_r}{P^2} e_1+\left( -\frac{\partial P}{P^2}+W \frac{P_r}{P^2}+\frac{W_r}{P}\right)e_4\\
%
%
%
 %
 %
 %
 %
 %
 %
 &[e_3, e_4]=\frac{P_r}{P^2}
 e_3+\left( -\frac{\partial_0 P}{P^2}+H \frac{P_r}{P^2}+\frac{H_r}{P} \right)e_4,
\end{align*}
where the subscripts 1, 2, 3  in the above expressions denote derivation with respect to the corresponding frame field \eqref{frame}. For example, $H_1$
means $\frac{1}{P} (\partial H -WH_r) $.

Now by using these commutator relations and Cartan's structure equations  
\begin{align*}
d\theta^i+\Gamma_k^i \wedge \theta^k=0,
\end{align*}
 for the metric \eqref{metric2} we find the connection forms listed  below:
\begin{align}
&\Gamma^1_4=(\frac{\mathrm{i}}{2P}+c^1_{\,14})\theta^1+\frac{1}{2}(c^3_{\,23}+c^4_{\,24})\theta^3,\label{gamma14}\\
&\Gamma^1_1=-c^2_{\,12}\theta^1-c^1_{\,12}\theta^2+\frac{1}{2}(c^2_{\,23}-c^1_{\,13}-c^4_{\,12})\theta^3+\frac{\mathrm{i}}{2P}\theta^4\label{gamma11}\\
&\Gamma^4_4=c^4_{\,34}\theta^3+c^3_{\,34}\theta^4-\frac{1}{2}c^3_{\,23}\theta^2-\frac{1}{2}c^3_{\,13}\theta^1+\frac{1}{2}c^4_{\,14}\theta^1+\frac{1}{2}c^4_{\,24}\theta^2\label{gamma44}\\
&\Gamma^3_1=(\frac{\mathrm{i}}{2P}-c^2_{\,24})\theta^2-\frac{1}{2}(c^3_{\,13}+c^4_{\,14})\theta^3\nonumber\\
%
%
%
%
%
%
%
%
%
&\Gamma^4_1=-c^2_{\,13}\theta^1
-\frac{1}{2}(c^4_{\,12}+c^2_{\,23}+c^1_{\,13})\theta^2
-c^4_{\,13}\theta^3
-\frac{1}{2}(c^4_{\,14}+c^3_{\,13})\theta^4\nonumber\\
%
%
%
%
%
%
&\Gamma^1_3=\frac{1}{2}(-c^4_{\,12}+c^1_{\,13}+c^2_{\,23})\theta^1+c^1_{\,23}\theta^2+c^4_{\,23}\theta^3+\frac{1}{2}(c^4_{\,24}+c^3_{\,23})\theta^4\nonumber
\end{align}
where $c^k_{\,\,mn}$ are the structure constants defined by $[e_m, e_n]=c^k_{\,\,mn}e_k$. 
 We also notice that $dg_{ij}=\Gamma_{ij}+\Gamma_{ji}$ and hence 
 $$\Gamma^1_2=\Gamma_{11}=0 \text{ and } \Gamma^3_4=\Gamma_{33}=0.$$

\begin{remark}
Note that, because of the choice of the coframe \eqref{coframe} the complex conjugate interchanges the indices $1$ and $2$ and keeps the indices $3$, $4$ unchanged, for example, ${\bar{\Gamma}}^1_4=\Gamma^2_4$.
\end{remark}

In the proposition below we compute the Ricci components of the quasi-Fefferman metric.
\begin{proposition}\label{ricpro}
Let $g$ be a quasi-Fefferman metric \eqref{metric1} on $\mathcal{M}^{\frac32}$ associated with a CR-manifold $M$ that admits a non-constant CR-function $\zeta$. Let $(\mu=d\zeta,\lambda)$ be a coframe for $M$ and $R_{ik}$ the components of the Ricci curvature with respect to an adapted frame. Then  
\begin{itemize}
\item[(i.)] $R_{44}=0$ is equivalent to 
\begin{equation}\label{solup}
P=\frac{a}{\cos( \frac{r+s}{2})}
 \end{equation}
where $a, s$ are arbitrary $r$-independent real functions.
\item[(ii.)] $R_{24}= 0$ is equivalent to 
 \begin{align}\label{sec}
\partial \log a^2+\I\partial s- 2x\e^{\I s}= -\frac{2 c}{3},
\end{align}
\item[(iii.)] $R_{22}= 0$ if and only if the equation 
 \begin{equation}\label{equ}
 \partial t+ t(c-t)=0
 \end{equation}
  is satisfied  where
   \begin{align}\label{t}
 t&=c + \partial \log  a^2-x\e^{\I s}.
 \end{align}
\end{itemize}
For an alternate coframe $(\mu',\lambda')$ the function $t$ changes to
 \begin{equation}\label{equt}
t'=\e^{-\tau-\I \theta}(t-\I \bar{h}).
\end{equation}
\end{proposition}
{\bf Proof.} 
To verify that the condition $R_{44}=0$ is equivalent to the function $P$ having the form \eqref{solup} we first notice that $R_{44}=2R^1_{414}$. We now consider the Cartan's  structure equation for the 1-form $\Gamma_4^1$
\begin{align*}
d\,\Gamma_4^1+\Gamma_k^1 \wedge \Gamma_4^k=R^1_{4k\ell}\theta^k \wedge \theta^\ell,\quad k<\ell,
\end{align*}
which takes the form 
\begin{align}\label{carp}
d\,\Gamma_4^1+\Gamma_1^1\wedge \Gamma_4^1+\Gamma_4^1 \wedge \Gamma_4^4=R^1_{4k\ell}\theta^k \wedge \theta^\ell,
\end{align} 
since $\Gamma_2^1=\Gamma_4^3=0$.
Substituting the 1-forms 
$\Gamma_4^1, \Gamma_1^1, \Gamma_4^4$ given by \eqref{gamma14},\eqref{gamma11},\eqref{gamma44} into \eqref{carp} and inspecting the coefficient of the 2-form $\theta^1 \wedge \theta^4$ in \eqref{carp}, we get the differential equation 
\begin{equation}\label{equp}
-4PP_{rr}+8P^2_r+P^2=0.
\end{equation}
The general solution of the differential equation \eqref{equp} has the form \eqref{solup}.

A similar argument shows that $R_{24}=0$ is equivalent to \eqref{sec}. 

To show that $R_{22}=0$ is equivalent to \eqref{equ}, we first notice that, due to $d\mu=0$, the structure functions $\alpha, \beta$ vanish. 
Then, by similar arguments as above,  $R_{22}=0$ becomes equivalent to the differential equation \eqref{equ}.

To verify  \eqref{equt}, we first  notice that
\begin{align*}
P'=\e^{-\tau}P=\frac{a\e^{-\tau}}{\cos(\frac{r+s}{2})}=
\frac{a'}{\cos(\frac{r'+s'}{2})}
\end{align*}
for all $r$ and $r'=r-\frac23\theta$. It follows $a'=\e^{-\tau} a$ and  $s'=\frac{2}{3}\theta+s$. Therefore, 
\begin{align*}
t'&=c'+\partial' \log a'^2-x'\e^{\I s'}=\e^{-\tau-\I \theta}\left(c-2\I \bar{h}+\partial (\tau+\I \theta) \right)+\e^{-\tau-\I \theta}(\partial \log a^2-2\partial \tau)-\e^{-\tau-\frac{5\I}{3}\theta} x\e^{\I s+\frac{2\I}{3}\theta}\\
&=\e^{-\tau-\I \theta}\left(c-2\I \bar{h}+\partial (\tau+\I \theta)+\partial \log a^2-\partial (\tau+\I \theta)+\I \bar{h}-x \e^{\I s} \right)=\e^{-\tau-\I \theta}\left(t- \I\bar{h}\right).\qquad \Box
\end{align*}

Another important ingredient is Jacobowitz's theorem \eqref{cls}, which uses the notion of the canonical bundle. 
\begin{definition}\label{defcan}
Let $(M,H,J)$ be a CR structure and $H^{\mathbb C}=H^{1,0}\oplus H^{0,1}$ the eigenspace decomposition of $J$. The  canonical bundle $\mathcal{K}$ over $M$ is the complex line bundle of complex-valued 2-forms with kernel $H^{0,1}$, i.e.
 $$\mathcal{K}=\{ \phi\in \Lambda^2(M)\otimes \mathbb C \mid  \phi(X,\cdot)=0 \quad \forall X\in H^{0,1}\}.$$
\end{definition}

If the CR-structure $(M,H,J)$ is given by the coframe $(\lambda, \mu)$ then the canonical bundle is spanned by  $\lambda \wedge \mu$. Using this representation the existence of a closed non-vanishing section can be reformulated as a $\bar{\partial}$-problem.

\begin{proposition}\label{clsecprop}
A CR-structure $(M,\lambda,\mu)$ admits locally a non-vanishing section of the canonical bundle if and only if the $\bar{\partial}$-problem
\begin{equation}\label{clsec}
\bar{\partial} \log \psi = -\bar{c},
\end{equation}
has a solution. Here $c$ is the structure function from \eqref{cofcon}
\end{proposition}

{\bf Proof.} Taking into account \eqref{cofcon},
$$d(\psi\, \mu\wedge \lambda)= \bar{\partial}\psi\, \bar{\mu}\wedge\mu\wedge \lambda-\bar{c}\psi \mu \wedge \bar{\mu}\wedge \lambda=(\bar{\partial}\psi+\bar{c}\psi)\,\bar{\mu}\wedge\mu\wedge \lambda$$
vanishes if and only iff \eqref{clsec} is satisfied with non-vanishing $\psi$. \hfill $\Box$ \medskip

While we do not assume the existence of a non-vanishing closed section of the canonical bundle a priori we show that this is a consequence of our assumptions.

We are now ready to prove our main theorem.

\begin{theorem}
A 3-dimensional CR-structure $(M,H,J)$ is (locally) embeddable if and only if  there exists an associated circle bundle $\mathcal{M}^{\frac32}$ with a quasi-Fefferman metric $g$
whose complexified Ricci tensor vanishes on the distribution of $\alpha$-planes. 
\end{theorem}
{\bf{Proof.}}
Let $M$ be a  CR-structure with coframe $(\mu,\lambda)$ and let $g$ be a quasi-Fefferman  metric  defined by \eqref{metric1} on   $\mathcal M^{\frac{3}{2}} $ for which $R_{22}= R_{24}= R_{44}= 0$. 
We consider the connection 1-form 
\begin{equation}\label{gam24}
\Gamma_{24}=\Gamma^1_{4}=\sigma \theta^1+\rho \theta^3
\end{equation}
from  \eqref{gamma14} where 
$$\sigma=\frac{\I}{2P}+\frac{P_r}{P^2}, \qquad \rho=-\frac{\bar{\partial}P}{P^2}+\frac{\overline{W}P_r}{P^2}-\frac{\bar{c}}{2P}+\frac{\overline{W}_r}{P}$$
Clearly, $\sigma \neq 0$ and therefore, 
the form $\Gamma_{24} \neq 0$.
Moreover, 
\begin{equation*}
\Gamma_{24} \wedge \bar{\Gamma}_{24}\neq 0,
\end{equation*}
since 
$$\Gamma_{24} \wedge \bar{\Gamma}_{24}= |\sigma|^2\,\theta^1 \wedge \theta^2 \quad \text{mod}\, \theta^3.$$
On the other hand, the conditions of the Goldberg-Sachs Theorem \eqref{goldsachs} with respect to the  \SF vector field $\partial_r$  are satisfied and, therefore  
$$\Psi_0= 
C_{4141}=R_{1414}=0
 \qquad \Psi_1=
C_{4341}
=\frac{1}{2}(R_{4341}+R_{1421})=0.$$
It follows
$$C_{4242}=\overline{C_{4141}}=0,\quad C_{4342}=\overline{C_{4341}}=0 $$
and furthermore, using the symmetries of the Riemann curvature,
 $$R_{ijk\ell}=R_{k\ell ij}, \qquad  R_{ijk\ell}=-R_{jik\ell}=-R_{ij\ell k},$$
 we get
\begin{equation}\label{golsa}
R_{2424}=0,\qquad R_{2434}+R_{2412}=0.
\end{equation}
Since
$$R_{44}=2R_{2414}, \quad R_{22}=2R_{2423}, \quad R_{24}=R_{2412}-R_{2434},$$
where $R_{ij}=R^k_{ikj}$ and $R_{ijk\ell}=g_{im}R^m_{jk\ell}$,
this shows that  the conditions $R_{44}=R_{22}=R_{24}=0$ are  equivalent to 
\begin{equation}\label{alpl}
R_{2414}=R_{2423}=R_{2412}-R_{2434}=0.
\end{equation}
Combining \eqref{golsa} and \eqref{alpl} yields
\begin{equation}\label{riecu}
R_{2412}=R_{2424}=R_{2414}=R_{2423}=R_{2434}=0.
\end{equation}
Therefore Cartan's structure equation \eqref{carp} for the connection 1-form $\Gamma_{24}=\Gamma^1_4$ becomes 
$$d\Gamma_{24}-(\Gamma_{12}  +\Gamma_{34}) \wedge\Gamma_{24}  =R_{24k\ell}\theta^k \wedge \theta^\ell=R_{2413}\theta^1 \wedge \theta^3.$$
Wedging the equation above with $\Gamma_{24}$ and taking into account that $\Gamma_{24}$ is a linear combination of $\theta^1$ and $\theta^3$, given by \eqref{gam24},
we conclude that 
$$d\Gamma_{24} \wedge \Gamma_{24}=0.$$
Now we can apply Lemma \eqref{firstcr} for the 1-form $\Gamma_{24}$ and deduce that locally
there exist complex functions $h\neq 0, \zeta$ 
such that 
$$\Gamma_{24}=h\,d\zeta \quad  \text{with} \quad d\zeta \wedge d\bar{\zeta}\neq 0.$$
Wedging the equation
$$h\,d\zeta=\Gamma_{24}=P(\sigma\mu+\rho\lambda)$$ 
by $\lambda \wedge \mu$ shows that
$$d\zeta \wedge \lambda \wedge  \mu=0.$$ 
Restricting the function $\zeta$ to the CR-manifold $M$, considered as a section $\{r=0\}$ of $\mathcal M^{\frac32}$, gives a CR-function there. 

Now we may assume that $\mu=d \zeta$. Since vanishing of the Ricci tensor on the $\alpha$-planes does not depend on the choice of an adapted frame, the conditions $R_{44}=R_{24}=R_{22}=0$ are still satisfied.

 If, coincidentally,   $t$ defined by \eqref{t} vanishes everywhere then from the equation \eqref{sec} it follows that 
 \begin{align*}
 \frac{4}{3} c=- \partial \log  a^2 +\I\partial s\\\intertext{and hence,} \partial \log (a^{\frac{3}{2}}\e^{-\frac{3}{4}\I s})=-c.
 \end{align*}
 Therefore, equation \eqref{clsec} has a solution $\psi= a^{\frac{3}{2}}\e^{\frac{3}{4}\I s}$
and, by Proposition \ref{clsecprop}, 
 the canonical bundle has a non-zero closed section. Now, by Theorem \eqref{cls},  the CR-structure is embeddable.  
 
Otherwise,  if  $t$ is not identically $0$, we
replace the complex coframe 1-form $\mu$ by another exact form $\mu'$. Consider
  \begin{equation}\label{dephi}
\varphi=\mu+\I \bar{t}\lambda. 
\end{equation}
Since $R_{22}=0$ we have
$$d\varphi \wedge \varphi=\I(\bar{\partial} \bar{t}+ \bar{t}(\bar{c}-\bar{t}\,))\mu \wedge \bar{\mu} \wedge \lambda=0.$$

Also $\varphi \wedge \bar{\varphi} \neq 0$ holds because 
$$\varphi \wedge \bar{\varphi}=\mu \wedge \bar{\mu}-\I t \mu \wedge \lambda-\I \bar{t}\bar{\mu}\wedge \lambda.$$
Thus the 1-form $\varphi$ satisfies the conditions of the Lemma \eqref{firstcr}. Consequently, there exist  complex-valued functions $b\neq 0, \eta$ such that 
\begin{equation}\label{crf}
\varphi=\mu+\I \bar{t}\lambda=b\, d \eta.
\end{equation}
Clearly,
$$d\eta \wedge d\bar{\eta}=\frac{1}{|b|^2}\varphi \wedge \bar{\varphi}\neq 0.$$
It follows from the definition of $\eta$ and $\varphi$ that $d\eta$ is a linear combination of $\mu$ and $\lambda$ and hence
$$d\eta \wedge \lambda \wedge \mu=0,$$
that is, $\eta$ is a CR-function.
Now we switch to the coframe $(\mu'=d\eta,\lambda')$  for which, because of \eqref{equt}, $t' \equiv 0$ everywhere. 
This reduces the second case to the first case and proves embeddability of $M$.

For the proof of the converse statement we assume that the CR-structure $M$ with adapted coframe $(\mu=d\zeta,\lambda)$ 
is embeddable. Then  the canonical bundle contains a non-zero closed section, i.e. there exists a non-zero complex function $\psi$ such that
$$\partial \log \bar{\psi}=-c.$$
We define real functions $a, s$ and a complex function $x$ as follows
$$\log a^2=\frac{4}{3} \text{Re}(\log \bar{\psi}), \qquad s=-\frac{4}{3}\text{Im}(\log \bar{\psi}), \qquad x=\e^{-\I s}(c+\partial\log a^2).$$
The metric  defined by
$$g=2P^2\big [\mu \bar{\mu}+\lambda (dr + W \mu +\overline{W}\bar{\mu} + H \lambda) \big ]$$
where 
$$P=\frac{a}{\cos( \frac{r+s}{2})}, \qquad W=\I x\e^{-\I r}-\frac{\I}{3}c, \qquad 
 $$
  and $H$ is any real function defined on $\mathcal M^{\frac{3}{2}}$,  is a quasi-Fefferman metric for $(M,\mu,\lambda)$ and, due to Proposition \eqref{ricpro},  $R_{44}=R_{24}=R_{22}=0$ is satisfied. \hfill $\Box$

\end{document}